\documentclass[12pt]{article}
\usepackage[all]{xy}
\usepackage{amsmath}
\usepackage{amsfonts}
\usepackage{amssymb}
\usepackage{amscd}
\usepackage{amsthm}
\usepackage{latexsym}
\usepackage{amsbsy}

\newtheorem{lem}{Lemma}[section]
\newtheorem{prop}{Proposition}[section]
\newtheorem{theorem}{Theorem}[section]

\newtheorem{definition}{Definition}[section]
\newtheorem{example}{Example}[section]

\begin{document}
\date{}
\title
{On the Generalized Cyclic Eilenberg-Zilber Theorem}
\author {M. Khalkhali,~~~ B. Rangipour,
\\ \footnotesize\text{~masoud@uwo.ca ~~~~brangipo@uwo.ca}
\\\footnotesize{Department of Mathematics }
\\\footnotesize{University of Western Ontario }}
\maketitle
\date{}
\begin{abstract}
We use the homological perturbation lemma to give an algebraic
 proof of the cyclic Eilenberg-Zilber theorem  for cylindrical modules.
\end{abstract} 

\section{Introduction}
The original Eilenberg-Zilber theorem (see \cite{jardin} for a recent account) 
states that if $X$ and $Y$ are simplicial abelian groups then the total complex
 of the bicomplex $X\otimes Y$ is chain homotopy equivalent  to the diagonal 
 complex of $X\otimes Y$. This result was then generalized by Dold and Puppe  
to bisimplicial abelian groups \cite{dp,jardin}: if $X$ is a bisimplicial abelian 
group then the total complex and the diagonal complex of $X$ are chain homotopy 
equivalent. This extension is important because in many examples (e.g. the 
bisimplicial group associated to a group action through its translation category),  
$X$  is not decomposable as the tensor product of two simplicial groups.

Thanks to the work of Connes \cite{ac},  one knows that a general setup for defining and studying  
cyclic homology is through cyclic modules. In order to define products and coproducts 
 in cyclic homology and to prove  K\"{u}nneth type formulas, 
 several authors, including Kassel \cite{k1},  Hood-Jones  \cite{hj},  and Loday  \cite{ld},
  have proved Eilenberg-Zilber type theorems for tensor products of cyclic modules 
  associated to algebras. A good reference that compares various methods 
  used by these authors is A. Bauval's article \cite{ba}.  
     In (\cite{ld}, page 130) one
 can find an Eilenberg-Zilber theorem for tensor products of two cyclic modules. 
  The most general result in this direction is stated by 
 Getzler and Jones in \cite{gj}. The proof is however topological in nature, and is 
  based on the method of acyclic   models.  In view of the importance of this result 
  for cyclic homology theory, for example in deriving a spectral sequence for the cyclic homology 
  of the crossed product algebra for the action of a group, as in  \cite{gj}, or for the action of a Hopf algebra,
  we felt it desirable  to give a purely algebraic proof of this fact.

In attempting to extend the proof  in \cite{ld} to this more general setup,
 we realized that  the cyclic shuffle map of \cite{ld} has no immediate extension to cylindrical 
(or even bicyclic) modules, while the shuffle and Alexander-Whitney maps have
 more or less obvious  extensions (see the Remark in Section  \ref{sec} for more on this). 
 It seems plausible that one should use a 
 different definition for cyclic shuffles. Instead,  we use the homological 
 perturbation lemma to obtain an algebraic proof for the cyclic Eilenberg-Zilber theorem for cylindrical modules.
 Our approach  is motivated by A. Bauval's work \cite{ba}, where a  perturbation lemma is used to give an alternative
  proof of the cyclic Eilenberg-Zilber theorem of  \cite{ld}. In using the perturbation lemma, one has to overcome
   two difficulties:  first,  showing that the first term in the perturbation formula is  identical with the boundary operator 
$B_t$ (Proposition  \ref{prop1}), and secondly,   proving that all higher order terms in the perturbation formula actually vanish 
(Theorem \ref{th}). We find it remarkable that these results continue to be true in  our cylindrical  module context. 
By making use of explicit formulas for the contracting homotopy in the generalized Eilenberg-Zilber theorem for 
bisimplicial 
modules, one can in principle, find an explicit formula for a generalized cyclic shuffle map.

We would like to thank  Rick Jardine and Jean-Louis Loday for informative  
discussions  on the subject of this paper.  We would like also to thank a referee 
whose critical  comments and suggestions led to a  substantial improvement in our 
presentation. 
     
\section{Preliminaries}
Let $k$ be a commutative unital ring. Recall that \cite{ft} a $\Lambda_{\infty}$-module
is a simplicial $k$-module $M=(M_n)_{n\geq 0}$ endowed, for each $n\geq 0$, with 
automorphisms   $\tau_n:M_n\longrightarrow M_n$, such that the following relations hold   
\begin{eqnarray*}
\delta_i\tau_n& =&\tau_{n-1}\delta_{i-1}  ,\;\;\;1\leq i\leq n,\\
\delta_0\tau_n&=&\delta_n,\\
\sigma_i\tau_n &=&\tau_{n+1}\sigma_{i-1},\;\;\;1\leq i\leq n,\\
\sigma_0\tau_n&=&\tau_{n+1}^2\sigma_n.
\end{eqnarray*}
Here $\delta_i$ and $\sigma_i$ are the faces and degeneracies of $M$. In case $\tau^{n+1}_n=1$ for all $n\geq 0$, we
 say that $M$ is a cyclic $k$-module. We denote the categories of $\Lambda_\infty$ (resp. cyclic) $k$-modules by
  $k\Lambda_\infty$ (resp. $k\Lambda$).

 For example, to each unital $k$-algebra $A$ and an algebra automorphism $g\in \mathcal{A}ut(A)$ one can associate a
  $\Lambda_{\infty}$-module $A^{\natural}_g$  with $A^{\natural}_{g,n}=A^{\otimes (n+1)}$,  and with faces, degeneracies
   and $\tau_n$  defined by 
\begin{eqnarray*}
\delta_i(a_0\otimes\dots \otimes a_n)&=&a_0\otimes a_1\otimes \dots \otimes a_ia_{i+1}  \dots \otimes a_n,\;\;\;0\leq
 i\leq n-1,\\
\delta_{n}(a_0\otimes\dots \otimes a_n)&=&(ga_n)a_0\otimes a_1\otimes \dots \otimes a_{n-1},\\
\sigma_i(a_0\otimes \dots\otimes a_n)&=& a_0\otimes \dots\otimes a_i\otimes 1\dots   \otimes a_{n},\;\;\;0\leq i\leq
 n,\\
\tau_n(a_0\otimes \dots\otimes a_n)&=& ga_n\otimes a_0 \dots\otimes a_{n-1}.
\end{eqnarray*}

The cyclic homology groups of  a cyclic $k$-module ${M}$ can be defined, among 
other ways, via the  bicomplex $\mathcal{B}(M)$ defined by 
\begin{eqnarray*}
\mathcal{B}_{m,n}(M)&=&M_{n-m} \;\;\text{if}\;\; n\geq m\ge 0, \\
\mathcal{B}_{m,n}(M)&=&0 \;\; \;\;\;\;\;\;\;\text{otherwise.} 
\end{eqnarray*}
The  vertical boundary operator is defined by
$$b=\sum_{i=0}^m(-1)^i\delta_i:M_m\rightarrow M_{m-1},\; $$
and the horizontal boundary operator is given by  
$$B=N \sigma_{-1} (1-t):M_m\rightarrow M_{m+1},$$
where  
$$t=(-1)^m\tau_m, \; \sigma_{-1}=\tau_m \sigma_m:M_m\rightarrow M_m, 
~\text{and}~ N=1+t+\dots +t^m .$$ One can check that $b^2=B^2=bB+Bb=0$.  
 The cyclic homology of $M$, denoted by  $HC_{\ast}(M)$, is defined to be 
the total  homology of the first quadrant  bicomplex  $\mathcal{B}(M)$. 

If $M$ is only a $\Lambda_{\infty}$-module, we can still define the operator $B$ 
as above and $B^2=0$ \cite{gj}, however $bB+Bb$ need not be zero. As in \cite{gj},  let $T=1-bB-Bb$.  

Recall that  a mixed complex $(C,b,B)$ is a chain complex $(C,b)$ with a  map of
 degree $+1$, $B:C_\ast\longrightarrow C_{\ast+1} $, satisfying $b^2=B^2=bB+Bb =0$ 
  \cite{ld}. To any mixed complex $C$ one associates  a bicomplex $\mathcal{B}C$ in the first quadrant,
  defined by  $\mathcal{B}C_{n,m}= C_{m-n}$ if $m\geq n\ge 0$ and $0$ otherwise, 
  with horizontal boundary $B$ and vertical boundary $b$.
By definition, the cyclic homology of $C$ is  $HC_{\ast}(C)=H_{\ast}({Tot} (\mathcal{B}C))$,
 and its Hochschild homology is, $HH_{\ast}(C)=H_{\ast}(C,b)$.  As for  cyclic 
 homology of algebras,   here also,  we have a short exact sequence of complexes,
$$
\CD 
0 @> >> (C,b)@> >> {Tot} (\mathcal{B}C) @> S>> {Tot}(\mathcal{B}C)\lbrack 2 \rbrack @> >> 0, \\
\endCD   
$$
where $S$ is the quotient map obtained by factoring by the first column.

An $S$-morphism of mixed complexes $f: (C,b,B)\longrightarrow (C',b',B')$ is a morphism of complexes
  $f: {Tot}(\mathcal{B}C) \longrightarrow {Tot}(\mathcal{B}C') $, such that  $f$ commutes with
   $S$. One can write  an  $S$-morphism as a matrix of maps
\begin{displaymath}
f=
\left( \begin{array}{cccc}
f^0&f^1&\ldots &\ldots \\
f^{-1}&f^0&f^1&\ldots \\
\ldots&f^{-1}&f^0&\ldots \\  
\vdots & \vdots & \vdots & \ddots
\end{array}\right)
\end {displaymath}
where $f^i :C_{\ast-2i}\longrightarrow C_{\ast}'$, and $\lbrack B,f^i\rbrack +\lbrack b,f^{i+1}\rbrack = 0$.
 
The benefit of  $S$-morphisms  may be seen in many cases when we do not have a 
cyclic map between  cyclic modules but we can have a $S$-morphism. Every $S$-morphism 
 induces a map $f_{\ast}: HC_{\ast}(C)\rightarrow HC_{\ast}(C')$ rendering the  following diagram commutative 
$$
\CD
0 @> >>C @> >> {Tot}(\mathcal{B}C)@> >> {Tot}(\mathcal{B}C)[2] @> >> 0  \\
   &  &       @Vf^0VV                 @VfVV                                 @Vf[2]VV    \\
0 @> >>C' @> >> {Tot}(\mathcal{B}C')@> >> {Tot}(\mathcal{B}C')[2] @> >> 0  \\
\endCD
$$
We have the following proposition which follows easily from the five lemma: 
\begin{prop}({\cite{ld}})
Let $f:(C,b,B)\longrightarrow (C',b',B')$ be an $S$-morphism of 
mixed complexes.Then $f_{\ast}^0:HH_{\ast} (C)\longrightarrow HH_{\ast}(C')$ 
is an isomorphism if and only if $f_{\ast}:HC_{\ast}(C)\longrightarrow HC_{\ast}(C')$ is an isomorphism.  
\end{prop}
Let $\Lambda_{\infty}$(resp.~$\Lambda$) be the $\infty$-cyclic (resp. cyclic)
 categories of Feigin-Tsygan \cite{ft}(resp. Connes \cite{ac}). We do not 
 need their actual definitions for this paper.  
\begin{definition} {(\cite{gj})} By a cylindrical $k$-module we mean a 
contravariant functor $X:\Lambda_{\infty} \times \Lambda_{\infty}  \rightarrow k$-mod,
  such that for all $p,q$, $\tau^{q+1}t^{p+1}=id: X_{p,q}\longrightarrow X_{p,q}$.
  More explicitly, we have a bigraded sequence of $k$-modules  $X_{n,m}$,  $n,m \geq 0 $,  
   with horizontal and vertical face, degeneracy and cyclic operators
\begin{quote}
$d_i: X_{n,m} \longrightarrow  X_{n-1,m}, \;\; \; \; n\geq i\geq 0 $,\\
$s_i: X_{n,m} \longrightarrow  X_{n+1,m}, \;\; \; \; n\geq i\geq 0 $,\\
$t: X_{n,m} \longrightarrow  X_{n,m}, $\\
  ${\delta}_i :X_{n,m} \longrightarrow  X_{n,m-1}, \;\; \; m\geq i\geq 0 $,\\
 $\sigma_i :X_{n,m} \longrightarrow  X_{n,m+1}, \;\; \; m\geq i\geq 0 $,\\
 $\tau :X_{n,m} \longrightarrow  X_{n,m}$, 
\end{quote}
 such that every vertical operator commutes with every horizontal operator and 
 vertical and horizontal operators satisfy the usual $\Lambda_\infty$-module  
 relations. Moreover, for all $p$ and $q$  the crucial relation
$\tau^{q+1}t^{p+1}=id : X_{p,q}\longrightarrow X_{p,q}$ holds. In this paper we denote the horizontal operators
 by $d_i$,$s_i$,$t$ and the vertical operators by $\delta_i$,$\sigma_i $,$\tau$.

\end{definition}

\begin{example}
Homotopy colimits of diagrams of simplicial sets are defined as the diagonal of
 a bisimplicial set (\cite{jardin}, page 199). We show that if the indexing category 
 is a cyclic groupoid and the functor has certain  extra properties, we can turn this 
 bisimplicial set into a cylindrical module.  
Let $I$ be a groupoid, i.e.,  a small category in which every morphism is an
 isomorphism. Recall from \cite{bu} that a cyclic structure $\varepsilon$ on 
 $I$ is a choice of morphisms $\varepsilon_i\in Hom(i,i)$ for all $i\in Obj(I)$,  
 such that for all $f:i\longrightarrow j$, $f\varepsilon_i=\varepsilon_jf$. Let
  $(I,\varepsilon)$ be a cyclic groupoid. We call a functor 
  $Z:I\rightarrow k\Lambda_{\infty}$ a cyclic functor if
   $Z(\varepsilon_i)\mid Z(i)_n=t_n^{n+1}\mid Z(i)_n$.
      To each cyclic functor $Z$ we associate a cylindrical  $k$-module  
$BE_IZ$, such that,
\[BE_IZ _{m,n}=\bigoplus_{i_0\overset{g_1}{\longrightarrow} i_1 \overset{g_2}{\longrightarrow}
 \dots \overset{g_m}{\longrightarrow} i_m} Z(i_m)_n.\]
 We   define the following cylindrical structure on $BE_IZ$:\\
\begin{eqnarray*}
 (\delta_j(x))_{i_0\overset{g_1}{\rightarrow} i_1 \overset{g_2}{\rightarrow} 
 \dots \overset{g_m}{\rightarrow} i_m}=
  \left \{ \begin{array}{lll}
(x)_{i_1\overset{g_2}{\longrightarrow} i_2 \overset{g_3}{\longrightarrow} 
\dots \overset{g_m}{\longrightarrow} i_m} & \textrm {if $ j=0 $} \\
(x)_{i_0\overset{g_1}{\longrightarrow} i_1 \overset{g_2}{\longrightarrow} 
\dots \overset{g_{j-1}}{\longrightarrow} 
i_{j-1}\overset{g_{j+1}\circ g_{j}}{\longrightarrow} i_{j+1}\overset{g_{j+2}}{\longrightarrow} \dots 
\overset{g_m}{\longrightarrow} i_m} & \textrm {if $ 1\leq j\leq m-1 $} \\
(g_m^{-1}(x))_{i_0\overset{g_1}{\longrightarrow} i_1 \overset{g_2}{\longrightarrow} \dots \overset{g_{m-1}}
{\longrightarrow} i_{m-1}} & \textrm {if $ j=m, $} \\
\end{array}\right.
\end{eqnarray*}\\
$$(\sigma_j(x))_{i_0\overset{g_1}{\longrightarrow} i_1 \overset{g_2}{\longrightarrow}
 \dots \overset{g_m}{\longrightarrow} i_m}
=(x)_{i_0\overset{g_1}{\longrightarrow} i_1 \overset{g_2}{\longrightarrow} \dots 
\overset{g_j}{\longrightarrow} i_{j}
\overset{id}{\longrightarrow} i_{j}\overset{g_{j+1}}{\longrightarrow} \dots
 \overset{g_m}{\longrightarrow} i_m},$$
$$(\tau(x))_{i_0\overset{g_1}{\longrightarrow} i_1 \overset{g_2}{\longrightarrow} 
\dots \overset{g_m}{\longrightarrow} i_m}
=(g_m^{-1}(x))_{i_m\overset{h}{\longrightarrow} i_0 \overset{g_1}{\longrightarrow}
 \dots \overset{g_{m-1}}{\longrightarrow} i_{m-1}},$$
where $h=(g_m\circ g_{m-1}\circ\dots \circ g_1)^{-1}$. 
The horizontal  structure is induced by the cyclic  structure of  $Z(i_m)$. One can check 
that $BE_IZ$ is a cylindrical module.

We apply the above construction to the following situation. Let $G$ be a (discrete) group acting by unital
 automorphisms on an unital $k$-algebra $A$. 
Let $I=G$ be the category with $G$ as its set of objects and  $Hom(g_1,g_2)=
\{h\in G\mid hg_1h^{-1}=g_2\}$. Define a
 cyclic structure $\varepsilon$ on $I$ by $\varepsilon_g=g$, for all $g\in G$.  Obviously $(I,\varepsilon)$ is a
  cyclic  groupoid. Define a functor $Z:I \rightarrow k\Lambda_\infty$ by $Z(g)=A^\natural_g$ and
 $Z(h):Z(g_1)\rightarrow Z(g_2)$ the map induced by $h$. It is a cyclic functor. It follows that    
  $BE_IZ$ is a cylindrical module. The cylindrical  module $X=BE_IZ$ can be identified as follows. We have   
$$X_{m,n}=\bigoplus_{G^{m+1}}A^{\otimes(n+1)}\cong kG^{\otimes (m+1)} \otimes A^{\otimes(n+1)}, $$
where $kG$ is the group algebra of the group $G$ over $k$.
The isomorphism is defined by 
\begin{center}
$\phi_{m,n}: BE_IZ _{m,n}\longrightarrow kG^{\otimes (m+1)} \otimes A^{\otimes(n+1)},$

$$\phi_{m,n}((a_0,a_1,\dots,a_n)_{i_0\overset{g_1}{\rightarrow}i_1\overset{g_2}
{\rightarrow}\dots \overset{g_m}{\rightarrow}i_m})=(i_m^{-1}g_m\dots 
g_1,g_1^{-1},\dots,g_m^{-1}\mid a_0,a_1,\dots,a_n).$$\\
\end{center}
 Under this isomorphism the   vertical and  horizontal  cyclic maps are given by:
\begin{multline*}
{\tau( g_0 , \dots , g_m \mid a_0 , \dots , a_n) }
= (g_0, \dots , g_m \mid g^{-1} \cdot a_n , a_0, \dots ,a_{n-1}), \\
\shoveleft{
\delta_i(g_0,\dots,g_m \mid a_0 , \dots , a_n)= (g_0,\dots,g_m \mid a_0,\dots , a_i a_{i+1},\dots , a_n), 
\;\;\; 0 \le i <n ,} \\
\shoveleft{\delta_n (g_0,\dots,g_m \mid a_0 , \dots , a_n)= (g_0,\dots,g_m \mid 
(g^{-1}\cdot a_n)a_0,\dots , a_{n-1}),} \\
\shoveleft{\sigma_i(g_0,\dots,g_m \mid a_0 , \dots , a_n) = (g_0,\dots,g_m\mid a_0,\dots , a_i , 1 ,
 a_{i+1},\dots , a_n),
\;\;\; 0 \le i \le n, }\\
\shoveleft{t( g_0 , \dots , g_m \mid a_0 , \dots , a_n)
=(g_m,g_0, \dots , g_{m-1} \mid g_m \cdot a_0, \dots
,g_m \cdot a_{n}),} \\
\shoveleft{d_i(g_0,\dots,g_m \mid a_0 , \dots , a_n)
=(g_0,\dots,g_i g_{i+1},\dots,g_m \mid a_0,\dots , a_n), \;\;\; 0 \le i <m ,} \\ 
\shoveleft{d_m (g_0,\dots,g_m \mid a_0 , \dots , a_n)
= (g_mg_0,g_1,\dots,g_{m-1} \mid g_m\cdot a_0,\dots , g_m\cdot a_{n}),}\\
\shoveleft{s_i (g_0,\dots,g_m \mid a_0 , \dots , a_n)
=(g_0,\dots,g_i,1,g_{i+1},\dots,g_m \mid a_0,\dots , a_n), \;\;\; 0 \le i \le m,} 
\end{multline*} 
where $g=g_0g_1\dots g_m.$
\end{example}
This shows that, in this particular case, our cylindrical module $BE_IZ$
 reduces to the cylindrical module associated in \cite{gj} to the action of a group on an algebra.   
\section{Proof of the main theorem}\label{sec}
Let $X $ be a cylindrical  $k$-module, and let $Tot(X)$
denote its total complex, defined by $Tot(X)_n=\oplus _{p+q=n}X_{p,q}$. 
 Consider the horizontal, and  vertical  $b$-differentials
  $b^h=\sum_{j=0}^n (-1)^jd_j$ and  $b^v=\sum_{i=0}^m (-1)^i\delta_i$,  and let
   $b=b^h+b^v$. Similarly we define $B^h$, $B^v$, and $B=T^vB^h+B^v$, 
   where $T^v=(1-b^vB^v-B^vb^v)$. The following lemma is proved in  \cite{gj}.

\begin{lem}
(Tot(X),b,B) is a mixed complex.
\end{lem}
We can  also define the diagonal $d(X)$ of a cylindrical  module $X$. It is a cyclic
 module with  $d(X)_n=X_{n,n}$,  and  the  cyclic  structure:  $ d_i\delta_i$ as  the $i$th
  face, $s_i\sigma_i$ as the $i$th degeneracy and  $t\tau$ as the cyclic map. Note that
   the cylindrical condition $\tau^{q+1}t^{p+1}=id: X_{p,q}\rightarrow X_{p,q}$ is needed in order to show
   that $(\tau t)^{n+1}=id: X_{n,n}\rightarrow X_{n,n}$. 
       Associated to this cyclic module we have a mixed complex $(d(X),b_d,B_d)$.

The following definition is from \cite{dp,we}. It extends the standard shuffle map,
 originally defined on the tensor product of simplicial modules,   to  bisimplicial modules.
\begin{definition}
Let $X$ be a bisimplicial module. Define 
$\nabla_{n,m} :X_{n,m}\longrightarrow X_{n+m,n+m}$ by 
\[\nabla_{n,m}=\sum_{\eta \in Sh_{m,n}}(-1)^{\eta}s_{\bar{\eta}(n+m)}\dots
 s_{\bar{\eta}(m+1)}\sigma_{\bar{\eta}(m)}\dots \sigma_{\bar{\eta}(1)},\]                 
where $Sh_{m,n}\subset S_{n+m}$, is the set of shuffles in the symmetric group
of order $n+m$, defined by $\eta \in Sh_{m,n}$ if and only if 
$\eta (1)< \eta (2) < \dots < \eta (m)$, and   $\eta (m+1)< \dots < \eta(m+n)$. 
Here $\bar{\eta}(j)=\eta (j)-1$, $1\leq j\leq n+m$. 
 We define the  shuffle map $$Sh :\bigoplus_{p+q=n}X_{p,q}\longrightarrow X_{n,n},$$ by
$$Sh= \sum _{p+q=n}\nabla _{p,q}.$$
\end{definition}
\begin{prop} 
The shuffle map $Sh:Tot(X)\rightarrow d(X)$ is a map of b-complexes of degree $0$.
\end{prop}
\begin{proof}
We must show that $ b_d\circ Sh =Sh\circ b^h+Sh\circ b^v $. 
All elements in the left hand side are of form 
$$d_i\delta_is_{\bar{\mu}(m+n)} \dots s_{\bar{\mu}(m+1)}\sigma_{\bar{\mu}(m)}\dots \sigma_{\bar{\mu}(1)}.$$
It would be better to divide these elements into five parts: 
\hspace{-5mm}
\begin{enumerate}
\item $i=0$, or $i=m+n$.
\item $1\leq i \leq m+n $, and   $i\in \{ \mu (1),\dots ,\mu (m)\}, i+1\in \{ \mu (m+1),\dots ,\mu (m+n)\}.  $ 
\item $1\leq i \leq m+n $, and   $i+1\in \{ \mu (1),\dots ,\mu (m)\},
 i\in \{ \mu (m+1),\dots ,\mu (m+n)\}.  $ 
\item  $1\leq i \leq m+n $, and   $i,i+1\in \{ \mu (1),\dots ,\mu (m)\}.  $ 
\item  $1\leq i \leq m+n $, and   $i,i+1\in \{ \mu (m+1),\dots ,\mu (m+n)\}.  $ 
\end{enumerate} 

For part $1$,  let $i=0$ (we leave to the  reader the rest of this case). We have
$$d_0\delta_0s_{\bar{\mu}(m+n)} \dots s_{\bar{\mu}(m+1)}\sigma_{\bar{\mu}(m)}\dots \sigma_{\bar{\mu}(1)}=
 s_{\bar{\bar{\mu}}(m+n)} \dots s_{\bar{\bar{\mu}}(m+1)}\sigma_{\bar{\bar{\mu}}(m)}\dots
  \sigma_{\bar{\bar{\mu}}(2)}d_0.$$
It is obvious that if we define 
$\rho (i)=\mu (i+1)-1$,  then $\rho$ is also a shuffle and the result is in $Sh\circ b^v $. 
\\ For case $2$, let $\mu(k)=i$, and $\mu(j)=i+1$, where $m+1\leq i\leq m+n$, and $1\leq j \leq m$.
  Now let $\alpha =\mu \circ (i,i+1)$.
           Then it is easy to check that $\alpha $ is also a shuffle and  
           we have $i+1\in \{ \alpha (1),\dots,  \alpha
            (m)\}$ and $i\in \{ \alpha (m+1),\dots ,\alpha (m+n)\}  $. On the other hand we have
             $$d_i\delta_is_{\bar{\alpha}(m+n)} \dots s_{\bar{\alpha}(m+1)}\sigma_{\bar{\alpha}(m)}\dots
              \sigma_{\bar{\alpha}(1)} = d_i\delta_is_{\bar{\mu}(m+n)} \dots
               s_{\bar{\mu}(m+1)}\sigma_{\bar{\mu}(m)}\dots \sigma_{\bar{\mu}(1)},$$
and $sign \mu =-sign\alpha $. So elements of case  $2$ cancel the elements of case $3$. 
Now let us do the case $4$. We assume $\mu(s)=i,\; \mu(s+1)=i+1$,  where $1\leq s\leq m$. We have 
$$ d_i\delta_is_{\bar{\mu}(m+n)}\dots s_{\bar{\mu}(s+2)}s_{i}s_{i-1}s_{\bar{\mu}(s-1)}
\dots s_{\bar{\mu}(m+1)}\sigma_{\bar{\mu}(m)}\dots \sigma_{\bar{\mu}(1)}$$
$$ = s_{\bar{\bar{\mu}}(m+n)}\dots s_{\bar{\bar{\mu}}(s+2)}s_{i-1}s_{\bar{\mu}(s-1)}
\dots s_{\bar{\mu}(m+1)}\sigma_{\bar{\theta}(m)}\dots \sigma_{\bar{\theta}(1)},$$
where
$$ 
\theta(j)=\left\{\begin{array}{ll}
\mu(j)\;\;\;\text{if}\;\;\;\mu(j)<i \\
\bar{\mu}(j)\;\;\;\text{if}\;\;\;\mu(j)>i+1.
\end{array}\right.
$$
 It is easy to show that the permutation 
$$
\small{
\left(\begin{array}{ccccccccccc}
1 \hspace{-2mm}&2\hspace{-2mm}&\hdots\hspace{-2mm} &m \hspace
{-2mm}&m+1\hspace{-2mm}&\hdots \hspace{-2mm}&s-1 \hspace{-2mm}&s 
\hspace{-2mm}& s+1\hspace{-2mm} & \hdots \hspace{-2mm}&m+n-1\\
\theta(1)\hspace{-2mm}&\theta(2)\hspace{-2mm}&\hdots \hspace{-2mm}&
\theta(m)  \hspace{-2mm}&\mu(m+1)  \hspace{-2mm}&\hdots \hspace{-2mm}&
\mu(s-1)\hspace{-2mm}&i-1 \hspace{-2mm}&\bar{\mu}(m+1)\hspace{-2mm}&\hdots
&\hspace{-2mm}\bar{\mu}(m+n)
\end{array}
\right)}
$$
is a $(m-1,n)$-shuffle. Similarly one can do case $5$ and then by counting   the proof is finished.
\end{proof}

In a similar way  to the shuffle map,  the  Alexander-Whitney map also extends to 
 bisimplicial modules \cite{we}. Define 
 $A_{p,q}:X_{n,n}\rightarrow X_{p,q}$, where $p+q=n$, by
 
  $$A_{p,q}=(-1)^{p+q}\delta_{p+1}\dots \delta_n \underbrace{ d _0\dots d _0}_{p\,  \text{times}},
  $$ 
  and let 
   $$A=\sum _{p+q=n}A_{p,q} :d(X)_n \longrightarrow Tot(X)_n.$$
   
Both maps $Sh$ and $A$  induce maps on the normalized complexes,  denoted by $\bar{A}: \bar{d}(X)\rightarrow
 \overline{Tot}(X)$ and $\overline{Sh}:\overline{Tot}(X)\rightarrow \bar{d}(X)$.\\ 
 \\
 {\bf Remark.}
   If $X=M\otimes N$ is the tensor product of two cyclic modules, one can define the cyclic shuffles \cite{ld}
 $$Sh': M_p\otimes N_q\longrightarrow M_{p+q+2}\otimes N_{p+q+2}.$$
 For example, if $M=A^\natural$ and $N=B^\natural$ are cyclic modules of associative unital algebras,
 we have $Sh'_{p,q}(x,y)=\sigma_{-1}(x)\perp \sigma_{-1}(y)$, where 
 $$\perp: A_p^\natural\otimes B_q^\natural \longrightarrow(A\otimes B)^\natural_{p+q},$$ 
 is defined by $$(a_0, a_1, \dots, a_p)\perp (b_0,  b_1,  \dots,  b_q)
 = \sum_\mu\mu^{-1}\cdot(a_0\otimes b_0,a_1\otimes 1, \dots a_p\otimes 1,1\otimes b_1,\dots 1\otimes b_q).$$  
 The summation is over all cyclic $(p,q)$-shuffles in $S_{p+q}$. The point is that one has to  use $\mu^{-1}$,
 as opposed to $\mu$ as  appears in \cite{ld} (page 127), in the above formula. This was also noticed by Bauval
  \cite{ba} and Loday (private communication). When $A$ and $B$ are commutative algebras, 
  the definition of cyclic shuffle,  and the fact that
 $(Sh, Sh'):\overline{Tot}(X)\longrightarrow\bar{d}(X)$ defines an $S$-map which is a quasi-isomorphism, are  
 essentially  due to G. Rinehart (cf. \cite{ld}).
 
 Now, one approach to prove the cyclic Eilenberg-Zilber theorem for cylindrical modules would be to extend 
 this cyclic shuffle map to the case where $X$ is not decomposable  to a tensor product. However, no
 natural extension  exists, exactly because one has to use $\mu^{-1}$  in the above formula for the cyclic shuffles.
  For example, with $p=2$ and $q=0$, $Sh'_{2,0}(a_0,a_1,a_2\mid b_0)$
  contains a term of the form $(1,a_0,a_2,1,a_1\mid 1,1,1,b_0,1)$, corresponding to the cyclic shuffle
    $\left(\begin{array}{cccc}
1 & 2 & 3 & 4 \\
1 & 3 & 4 & 2
\end{array}
\right)$. It is clear that  a term like $(a_0,a_2,a_1)$ can not be produced from $(a_0,a_1,a_2)$ 
via the cyclic or simplicial maps.  
\begin{prop}\label{prop}
$\overline{Sh}$ and $\bar{A}$ define a deformation retraction of $\bar{d}(X)$ to
 $\overline{Tot}(X)$, i.e. there is a homotopy $h:\overline{d}(X)\rightarrow \overline{d}(X)$ 
 such that $\overline{A}\circ\overline{Sh}=1$ and $\overline{Sh}\circ\overline{A}=1+bh+hb$.
\end{prop}
\begin{proof}
The existence of $h$ is part of the generalized Eilenberg-Zilber theorem \cite{jardin}. 
We just prove that  $\overline{A}\circ\overline{Sh} =1 $. Let us calculate the action 
of $\overline{A}\circ\overline{Sh}$  on a typical element $x\in X_{p,q}$. 
For every $\overline{A_{p',q'}}$ with $(p',q')\neq(p,q)$,  we have 
$\overline{A_{p',q'}}\circ\overline{Sh}(x)=0$, so we should only check the 
identity $\overline{A_{p,q}}\circ\overline{Sh}(x)=x$. For  $\mu\in Sh_{p,q}$ 
for simplicity let  $\mu .x=(-1)^\mu s_{\bar{\mu}(n)}\dots  s_{\bar{\mu}(p+1)}\sigma_{\bar{\mu}(p)}
\dots \sigma_{\bar{\mu}(1)}(x)  $, and $\mu_{p,q}=(q+1,\dots,n,1,\dots,q)$. 
Then it is not difficult to verify that $\overline{A_{p,q}}(\mu.x)=0$
 for every $\mu\neq\mu_{p,q}$, and $\overline{A_{p,q}}(\mu_{p,q}.x)=x$.
\end{proof}  
 The first applications of perturbation theory to cyclic homology go back to the work of C. Kassel in \cite{k}.
Our proof of the generalized Eilenberg-Zilber theorem, Theorem \ref{th} below, is also based on homological 
perturbation theory. We recall the necessary definitions and  results from \cite{ba,k}.  
 A chain complex $(L,b)$ is called a
 {\it deformation retract} of a chain complex $(M,b)$ if there are chain maps 
$$
\CD
  L@>g>>M @>f>> L 
\endCD 
$$
and a chain homotopy $h:M\rightarrow M$ such that
$$fg=id_L \;\;\; \text{and }\;\;\;g f=id_M+bh+hb. $$
 The retraction  is called {\it special} if in addition 
$$hg = fh=h^2=0.$$

It is easy to see that any retraction data as above can be replaced by a special retraction \cite{k}. 
Now we perturb the differential of  the \textquoteleft\textquoteleft bigger" complex  $M$ to $b+B$ so that $(b+B)^2=0$.
 It is natural to ask whether the differential of $L$ can be perturbed to  a new differential $b+B_\infty$ so that $(L,b+B_\infty)$ is a
  deformation retraction of $(M,b+B)$. The homological 
    perturbation lemma asserts that, under suitable conditions, this is possible. More precisely,
        assume the retraction is special, $L$ and $M$ have bounded below increasing
   filtrations, $g$ and  $f$  preserve the filtration and $h$ decreases the filtration. Then it is easy to check that
    the following formulas are well defined and define a special deformation retract of $(M,b+B)$ to $(L,b+B_\infty)$:
\begin{eqnarray*}
h_\infty &=&h\sum_{m\geq 0}(Bh)^m,\\
g_\infty &=&(1+h_\infty B)g,\\
f_\infty &=&f(1+Bh_\infty ),\\
B_\infty&=&f(1+Bh_\infty)Bg.
\end{eqnarray*}

To apply the perturbation lemma to our problem, we need to know that the perturbed
 differential $b+B_\infty$ coincides with the existing differential. This means 
 we have to show that the first term in the series for $B_\infty $
coincides with $B_t$ and all other terms vanish.   
The next proposition verifies the first part of
the claim. It  is a generalization of Lemma IV.1 in \cite{ba}.

\begin{prop}\label{prop1}
Let $X$ be a cylindrical module and let $B_t=T^vB^h+B^v$ and
  $B_d=B^hB^v$ be the total and diagonal $B$-differentials on 
  $\overline{Tot}(X)$ and $\overline{d}(X)$, respectively.  Then $\overline{A}B_d\overline{Sh}={B_t}. $
\end{prop}   
\begin{proof}
Let $x\in X_{p,q}$. As in  the proof of Proposition  \ref{prop},  we have $\overline{Sh}(x)=\sum_{\mu\in Sh_{p,q}}\mu .x.$

 Let $n=p+q$. In  
 $$\overline{A}B_d\overline{Sh}(x)=\sum_{r+s=n+1}\overline{A}_{r,s}B_d\overline{Sh}(x),$$
 because we are working with normalized chains, 
   all parts are  zero except for $r=p+1, s=q$ or $r=p, s=q+1$. 
  We denote the first part by $S_1$ and the second part  by $S_2$. We show
that  $S_1=T^vB^h(x)$ and $S_2=B^v(x)$.

For  $0\leq i\leq n$ and $\mu\in Sh_{p,q},$ let
 $$ S_1(i,\mu)=\overline{A}_{p+1,q}\tau_{n+1}\sigma_nt_{n+1}s_n\tau_n^it_n^i(\mu .x).$$
 The reader can easily check that  $S_1(i,\mu)=0$
 for all $0\leq i\leq {q-1}$ and all $\mu\in Sh_{p,q}$. For the rest of the 
 elements in $S_1$,  we have $S_1(i,\mu)=0$ for all $q\leq i\leq n$ and all $\mu\neq\mu_{p,q,i}$,  where \\
$$\mu_{p,q,i}=(1,2,\dots ,n-i,n+q-i+1,\dots ,n,n-i+2,\dots ,n+q-i).$$ 
 Now, we have   $S_1(i,\mu_{p,q,i})=(-1)^{(i-q)p}t_{p+1}s_pt_p^{i-q}\tau_q^{q+1}$. 
 We have shown that $S_1=T^vB^h(x).$ 
 
 Next, we compute $S_2$. For $0\le i\le n$  and $\mu\in Sh_{p,q}$ let 
  $$ S_2(i,\mu)=\overline{A}_{p,q+1}\tau_{n+1}\sigma_nt_{n+1}s_n\tau_n^it_n^i(\mu .x).$$

 For
$q+1\leq i\leq p+q$ and all $\mu \in Sh_{p,q}$ we have $S_2(i,\mu)=0$,  and 
if we denote $$\alpha_{i,p,q}=(q-i+1,\dots,n-i,1,2,\dots,q-i,n-i+1,\dots,n),$$ 
 then $S_2(i,\mu)=0$ for all $0\leq i\leq q$ and all $\mu\neq \alpha_{p,q,i}$. 
 Finally,  we have $S_2(i,\alpha_{p,q,i})=(-1)^{iq}\tau_{q+1}\sigma_q\tau_q^i$ for $0\leq i\leq q$. We have shown 
 that $S_2=B^v(x)$. The proposition is proved.

\end{proof}

Now we are in a position to combine the  perturbation lemma with the above proposition to prove: 
\begin{theorem}\label{th}
Let $X$ be a cylindrical module. Then there exists an $S$-map of mixed complexes, $f: Tot(X)\longrightarrow d(X)$,
such that $f_0=Sh$ is the shuffle map and $f$ is a quasi-isomorphism.
 
\end{theorem}
\begin{proof}
It suffices to prove the statement for the normalized complexes. By Proposition  \ref{prop},
 $(\overline{Tot}(X),b)$ is a deformation retract of $(\bar{d}(X),b)$. So, applying 
 the perturbation lemma, we have to show that all the extra terms in the perturbation 
 series vanish. Now the normalized homotopy operator is
 induced from the original homotopy operator $h: X_{n,n}\longrightarrow X_{n+1,n+1} $. Dold and Puppe 
show that the operator 
 $h$ is universal, in the sense that  it is a linear combination (with integral coefficients) of 
 simplicial morphisms of $X$ (page 213, Satz 2.9 in \cite{dp} ). One  
 knows that any order preserving map $\Phi :\lbrack n \rbrack \rightarrow \lbrack m \rbrack$ between finite 
 ordinals can be uniquely decomposed as,
 $$\Phi=\delta_{i_1}\delta_{i_2}\dots \delta_{i_r}\sigma_{j_1}\sigma_{j_2}\dots\sigma_{j_s},$$ 
 such that $i_1\le i_2\le\dots \le i_r$ and $j_1<j_2<\dots <j_s$ and $\delta_{i_k}$ are 
 cofaces  and $\sigma_{j_k}$ are codegeneracies
  (cf. e.g. Loday \cite{ld}, page 401). On dualizing this  and 
 combining it with the Dold-Puppe  result, it follows that the homotopy operator 
 $h:X_{n,n}\longrightarrow X_{n+1,n+1}$ is a linear combination of operators of the  
 form $\sigma_i^v\sigma_i^h\circ g$, where $\sigma_i^v$, $\sigma_i^h$ are vertical and horizontal 
 degeneracy operators and $g$ is another operator whose specific form is not important for the sake  
 of this argument. Now we look at the perturbation formula. We are claiming that 
 the induced operator on the normalized chains 
 $$\bar A\circ B_d\circ h_\infty \circ B_d\circ \overline{Sh}: \overline{Tot}(X)\rightarrow \overline{Tot}(X), $$  
 where $h_\infty =h\sum_{m\ge 0 } (B_dh)^m$ is zero.
 This follows from the above observation since 
 \begin{center}
 Image $(B_d\circ h_\infty\circ  B_d \circ Sh)\subset $ Image $(B_d\circ h)\subset$ degenerate chains.
   \end{center}
\end{proof}


\end{document}